\documentclass[11pt]{amsart}
\usepackage{color}
\usepackage{amsfonts, amsmath, amssymb}
\usepackage{wasysym}
\usepackage{graphicx}

\definecolor{RedClr}{rgb}{1,0,0}
\definecolor{BlueClr}{rgb}{0,0,1}
\definecolor{TextColor}{rgb}{0,0,0.5}
\definecolor{Violet}{rgb}{0.5,0,1}
\definecolor{Bordeaux}{rgb}{1,0.3,0.4}

\newtheorem{thm}{Theorem}[section]

\theoremstyle{definition}

\theoremstyle{remark}
\newtheorem{rem}[thm]{Remark}
\numberwithin{equation}{section} \theoremstyle{quest}

\numberwithin{equation}{section} \theoremstyle{prob}

\numberwithin{equation}{section} \theoremstyle{answer}

\numberwithin{equation}{section}

\begin{document}
		%
		\title{A Simple Differential Geometry for Networks and its Generalizations}
		%
		%
		\author{Emil Saucan$^*$, Areejit Samal$^{**}$ and J\"{u}rgen Jost$^{\dag,\ddag}$} 
		\address{$^*$Department of Applied Mathematics, ORT Braude College, Karmiel 2161002, Israel 
			          \and $^{**}$The Institute of Mathematical Sciences (IMSc), Homi Bhabha National Institute (HBNI), Chennai 600113, India
		              \and $^{\dag}$Max Planck Institute for Mathematics in the Sciences, Leipzig 04103, Germany 
	              \and $^{\ddag}$Santa Fe Institute; 1399 Hyde Park Road Santa Fe, New Mexico 87501, USA}
 
	              \email{semil@braude.ac.il, asamal@imsc.res.in, jost@mis.mpg.de}

	
	\maketitle
	
	\begin{abstract}
		Based on two classical notions of curvature for curves in general metric spaces, namely the Menger and Haantjes curvatures,  
		we introduce new definitions of sectional, Ricci and scalar curvature for networks and their higher dimensional counterparts. 
		These new types of curvature, that apply to weighted and unweighted, directed or undirected networks, are far more intuitive 
		and easier to compute, than other network curvatures. In particular, the proposed curvatures based on the interpretation of 
		Haantjes definition as geodesic curvature, and derived via a fitting discrete Gauss-Bonnet Theorem, are quite flexible. We also 
		propose even simpler and more intuitive substitutes of the Haantjes curvature, that allow for even faster and easier computations 
		in large-scale networks.
	\end{abstract}

	

\section{Introduction}
Complex networks are usually modeled as graphs. One can see such a
graph as a combinatorial object via its adjacency matrix that encodes
the neighborhood relations, or as a metric space, based on the
distance between vertices and the alternative paths between them.
For a long time, the combinatorial approach dominated, even though one
of the most classical and widely employed combinatorial indices, the
{\em clustering coefficient}, represents, in fact, a discretization of
the classical Gauss curvature \cite{EM}. Recently, 
the geometric approach has gained considerable momentum. This came
about because it was realized that notions that originated in
differential and Riemannian geometry possess a wider geometric meaning
that extends to metric spaces, and therefore, in particular, to
graphs. Notions of Ricci curvature have been particularly successful. 
In particular, Ollivier's Ricci curvature \cite{Ol1} has become an
established  method in the study of complex networks in their various
avatars \cite{NLGGS,Allen1,NLLG,S++,AGE}. We have also proposed 
\cite{SMJSS,WSJ17} yet another approach towards the introduction of 
Ricci curvature in the study of networks that originates with 
Forman's paper \cite{Fo}, and compared these two notions of Ricci 
curvature \cite{SMJSS,SSGLSJ}.

Ricci curvature always involves some averaging. On graphs, it is
assigned to edges. Ollivier's Ricci curvature is based on the concept of optimal
transport, but is prohibitively hard to compute for large networks and
for general weights. In contrast, Forman's version is extremely simple
to compute. However, it is less intuitive than Ollivier's one, as it
is based on a discretization of the so called
Bochner-Weizenb{\"o}ck-formula \cite{J}. In any case, the notion of
Ricci curvature is not the most elementary concept of geometry, and
therefore, the underlying geometric intuition may elude large parts of
the  active communities of engineers, social scientists, and
biologists interested in network analysis. 

Therefore, here we take a step back and start with the most elementary
notion of curvature, that of a curve. The idea of extending the
corresponding differential geometric notion to more general metric
spaces goes already back to Menger \cite{Me} and Haantjes \cite{Ha}.
By partially extending ideas that we already applied in the context of
Imaging and Graphics (and $PL$ manifolds in general) \cite{S15,Sa18,GS} 
we show that they allow us to naturally define geometrically intuitive 
notions of  curvatures for networks, multiplex-networks and 
hypernetworks, both weighted and unweighted, directed as well as 
undirected ones.

Menger's approach is the following. The curvature of a circle of
radius $R$ in the plane is $\frac{1}{R}$, and Menger then uses this 
to assign curvature values to other curves. On a graph, we can apply that
to triangles, and in order to get a Ricci type curvature for an edge,
we could simply average over all triangles containing that edge. 

Haantjes' approach is based on the observation that the curvature of
an arc depends on the difference between its length and the distance
between its endpoints. The larger that difference, the larger the
curvature. In an unweighted  graph, between two adjacent vertices
there is an edge, and this length is assigned length 1, and so this is
the distance between its endpoints. There may be alternative longer
paths between them, and they then get assigned a correspondingly
larger curvature. Haantjes curvature, in contrast to Menger, is not 
restricted to triangles. In this case, one can again average to get 
a Ricci type curvature. Haantjes cutrvature has two distinct advantages. 
Firstly, it is applicable to any 2-cells, not just to triangles. Secondly, 
it is better suited as discrete version of the classical geodesic curvature 
(of curves on smooth surfaces). As we shall see, it is therefore applicable 
to general networks, without {\it any assumption on the background geometry}. 
In fact, not only can it be employed in networks of {\em variable curvature}, 
it can be used to {\em define} the curvature of such a discrete space (see Section 3).

\section{Menger Curvature} \label{section:Menger}
The simplest, most elementary  manner of introducing curvature  in
metric spaces is due to Menger \cite{Me}. One simply defines  the
curvature  $K(T)$ of a triangle $T$  (metric triple of points) with
sides of lengths $a,b,c$  as
$1/R(T)^2$, where $R(T)$ is the radius of the circle circumscribed to
the triangle. An elementary computation yields 
\begin{equation} \label{eq:K-Euc}
\kappa_{M,E}(T) = \frac{1}{R(T)} = \frac{abc}{4\sqrt{p(p-a)(p-b)(p-c)}}\,,
\end{equation}
where $p = (a+b+c)/2$ denotes the half-perimeter.

However, there is conceptual problem with the above definition which 
utilizes the geometry of the Euclidean plane. In the general setting of 
networks, it is not natural to assume a Euclidean background. This is 
analogous to the geometry of surfaces, where the metric need not be 
Euclidean, but could be hyperbolic, spherical, or of varying Gauss 
curvature. For example, embedding networks in hyperbolic plane and 
space is becoming  quite common \cite{BR,Kr++,Z++}.

Of course, one may formulate a hyperbolic or spherical analogue of
\eqref{eq:K-Euc} (see, e.g. \cite{Ja}). The hyperbolic version is  
\begin{equation} \label{eq:K-Hyp}
\kappa_{M,H}(T) = \frac{1}{\tanh{R(T)}} =
\frac{\sqrt{\sinh{p}\sinh{(p-a)}\sinh{(p-b)}\sinh{(p-c)}}}{2\sinh{\frac{a}{2}}\sinh{\frac{b}{2}}\sinh{\frac{c}{2}}}\,;
\end{equation}
whereas the spherical one is
\begin{equation} \label{eq:K-Sph}
\kappa_{M,S}(T) = \frac{1}{\tan{R(T)}} = \frac{\sqrt{\sin{p}\sin{(p-a)}\sin{(p-b)}\sin{(p-c)}}}{2\sin{\frac{a}{2}}\sin{\frac{b}{2}}\sin{\frac{c}{2}}}\,.
\end{equation}
Note that, in the setting of networks, the constant factors ``4'' and  ``2'', 
respectively, appearing in the denominators of the formulas above are less relevant 
and they can be discarded in this context. 

Hyperbolic geometry is considered better suited to represent the background 
network geometry as it captures the qualitative aspects of networks of exponential 
growth such as the World Wide Web, and thus, it is used as the setting for variety 
of purposes. However, spherical geometry is usually not considered as a model 
geometry for networks because that geometry has finite diameter, hence finite growth.
However, spherical networks naturally arise in at least two instances. The first 
one is that of global communication, where the vertices represent relay stations, 
satellites, sensors or antennas that are distributed over the geo-sphere or over a 
thin spherical shell that can -- and usually is -- modeled as a sphere. The second
one is that of brain networks, where the cortex neurons are
envisioned, due to the spherical topology of the brain, as being
distributed on a sphere or, in some cases, again on a very thin (only
a few neurons deep) spherical shell, that can also be viewed as
essentially spherical.

One can also devise an analogous, although less explicit formula in
spaces of variable curvature, but in network analysis, it is not clear
where that background curvature should come from. After all, the
purpose here is to define curvature, and not take it as given.

As defined, the Menger curvature is always positive. This may not be
desirable, as in geometry, the distinction between positive and
negative curvature is important. For directed networks, however,  
a sign $\varepsilon(T) \in \{-1,0,+1\}$ is naturally attached to an 
directed triangle $T$ (see Figure 1), and the sectional-Menger curvature 
of the directed triangle is then defined, in a straightforward manner as
\begin{equation}
\kappa_{M,O}(T) = \varepsilon(T)\cdot\kappa_M(T),
\end{equation}
where $\kappa_M$ could be the Euclidean, the hyperbolic or the 
spherical version, accordingly to the given setting. 
%
\begin{figure} 
	\centering
	\includegraphics[height=4.3cm]{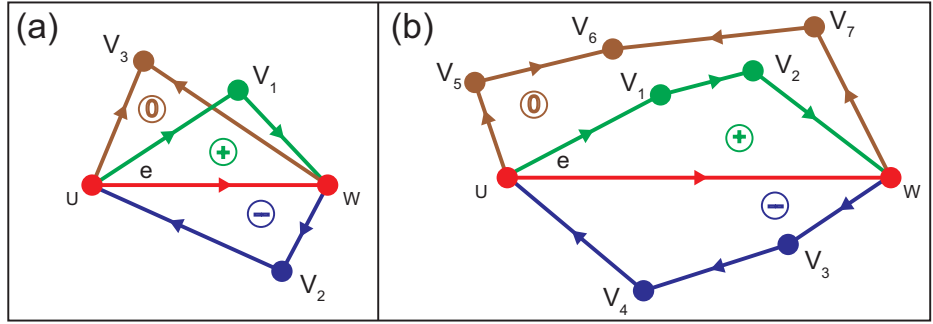}
	\caption{Sign convention for {\bf (a)} directed triangles and {\bf (b)} directed polygons.}
\end{figure}
We can then define  the directed Ricci curvature of an edge by
averaging  as in differential or piecewise linear geometry  as
\begin{equation}
\kappa_M(e) \stackrel{\rm not}{=} {\rm Ric}_M(e) = \sum_{T_e \sim e}\kappa_{M,O}(T_e)\,,
\end{equation}  
where $T_e \sim e$ denote the triangles adjacent to the edge $e$; and
\begin{equation} \label{eq:Ricci-Menger}
\kappa_M(v) \stackrel{\rm not}{=} {\rm scal}_M(v) = \sum_{e_k \sim v} {\rm Ric}_M(e_k) = \sum_{T \sim v}{\kappa_{M,O}(T)}\,;
\end{equation}

where,  $e_k \sim v, T \sim v$ stand for all the edges adjacent to the vertex 
$v$ and all the triangles $T$ having $v$ as a vertex, respectively. 
\begin{rem}
	${\rm Ric}_M(e)$ captures, in keeping with the intuition behind $\kappa_M(T)$ 
	the geodesic dispersion rate aspect of Ricci curvature. 
	(See \cite{SSGLSJ} for a succinct overview of the different aspects of Ricci curvature.)
\end{rem}

As already indicated, we see two drawbacks for Menger curvature as a
tool in network analysis. It depends on a background geometry model, 
and it naturally applies only to triangles, but not to more general
2-cells. Therefore, we next turn to Haantjes curvature.


\section{Haantjes curvature} \label{section:Haantjes}

Haantjes \cite{Ha} defined metric curvature  by comparing the ratio
between the length of an arc of curve and that of the chord it
subtends. More precisely, if $c$ is a curve in a metric space with metric 
$d$, and $p, q, r$ are points on $p$ between $q$ and $r$, the 
Haantjes curvature is defined as
\begin{equation}                         \label{eq:haantjes2}
\kappa_{H}^2(p) = 24\lim_{q,r \rightarrow
	p}\frac{l(\widehat{qr})-d(q,r)}{\big(d(q,r)\big)^3}\,\,;
\end{equation}
where $l(\widehat{qr})$ denotes the length -- in the intrinsic 
metric induced by $d$ -- of the arc $\widehat{qr}$.

In the network case, $\widehat{qr}$ is replaced by a path $\pi
= v_0,v_1,\ldots,v_n$, and the subtending chord by  $\bar{e} =
\overline{v_0v_n}$. Clearly, the limiting process has no meaning in
this discrete case. Furthermore, the normalizing constant 24 which
ensures that the limit will coincide, in the case of smooth planar
curves, with the classical notion, is superfluous in this
setting. This leads to the following definition of the 
{\it Haantjes curvature of a path} $\pi$: 
\begin{equation}                         \label{eq:haantjes-path}
\kappa_{H}^2(\pi) = \frac{l(\pi)-l(v_0v_n)}{l(v_0v_n)^3}\,\,;
\end{equation}
where, if the graph is a metric graph, $l(v_0v_n) = d(v_0v_n)$. 
In particular, in the case of the combinatorial metric, we obtain that, 
for $\pi$ as above, $\kappa_H(\pi) = \sqrt{n-1}$. 

Clearly, one can extend the above definition to directed paths in the 
same manner as done for the Menger curvature, namely 
\begin{equation}
\kappa_{H,O}(T) = \varepsilon(\pi)\cdot\kappa_H(T)\,;
\end{equation}
for every directed path $\pi$, where $\varepsilon \in \{-1,0,+1\}$ denotes 
the orientation of $\pi$. 

\subsection{A Local Gauss-Bonnet Theorem and the Curvature of 2-Cells}

Because of its advantages over Menger curvature, we shall now use
Haantjes curvature to define scalar and Ricci curvatures of networks. 
The basic idea here is to adapt  the {\it local Gauss-Bonnet Theorem} 
to this discrete setting.  Recall that, in the classical context of smooth 
surfaces, the theorem states that
\begin{equation}  \label{eq:SmoothGB}
\iint_DKdA + \sum_0^p\int_{v_i}^{v_{i+1}}k_gdl + \sum_{0}^p\varphi_i = 2\pi\chi(D)\,;
\end{equation}
where $D \simeq \mathbb{B}^2$ is a (simple) region in the surface, having as boundary $\partial D$ 
a piecewise-smooth curve $\pi$, of vertices (i.e. points where $\partial D$ is not smooth) 
$v_i, i = 1,\dots,n$, ($v_n = v_0$); $\varphi_i$ denotes the external angles of $\partial D$ 
at the vertex $v_i$; and $K$ and $k_g$ denote (as usually) the Gaussian and geodesic curvatures, 
respectively. 

We should first note that, in the absence of a background curvature, the very notion of angle 
is undefinable. Therefore, for abstract (non-embedded) cells, no ``honest'' notion of angle exists. 
Therefore, the last term on the left side of (\ref{eq:SmoothGB}) above has no proper meaning, 
thus should be discarded. Indeed, the distances between non-adjacent vertices on the same cycle 
(apart from the path metric) are not defined, thus the third term in the left side of formula 
(\ref{eq:SmoothGB}) vanishes. 

We first concentrate on the case of combinatorial graphs. For this type of networks, 
i.e. endowed with the combinatorial metric, the area of each cell is commonly taken 
as being equal to 1. Moreover, one assumes (quite naturally) that curvature is constant 
on each cell. Therefore, the first term in the left side of (\ref{eq:SmoothGB}) 
reduces simply to $K$. In addition, given that $D$ is a 2-cell, thus $\chi(D) = 1$. 
It follows, that in such a setting we obtain
\begin{equation}
K = 2\pi - \int_{\partial D}k_gdl\,.
\end{equation}
It is tempting to next consider $\partial D$ as being composed of
segments (on which $k_g$ vanishes), except at the vertices, thus
rendering the expression above as 
\begin{equation}\label{10}
K = 2\pi - \sum_{1}^n\kappa_H(v_i)\,.
\end{equation}
However, in general weighted graphs, one can not define a
(non-trivial) Haantjes curvature for each of the vertices since, as
already noted above, no proper distance  between the vertices
$v_{i-1}$ and $v_{i+1}$ can be implicitly assumed (apart from the one
given by the path metric, which would produce trivial 0 curvature at
$v_i$). In fact, in this general case, neither can the arc (path) $\pi
= v_0v_1\ldots v_n$ be truly viewed as smooth. Therefore,  we have no
choice but to replace the right term in \eqref{10} above by $\kappa_H(\pi)$, 
where it should be remembered that $\pi$ represents the path $v_0,v_1,\ldots,v_n$,  
of chord $\bar{e} = \overline{v_0v_n}$. 

We can now define the ({\it Haantjes}) {\it sectional curvature of a
	2-cell} $\mathfrak{c}$.  Given an edge $e = (u,v)$ and  a cell
$\mathfrak{c}$, $\partial \mathfrak{c} = (u=v_0,v_1,\ldots,v_n=v)$ 
({\it relative to the edge $e \in \partial \mathfrak{c}$}), we put 
\begin{equation} \label{eq:K=-kH}
K_{H,e}(\mathfrak{c}) = 2\pi -\kappa_{H,e}(\pi)\,;
\end{equation}
where $\pi$ denotes the path $v_0,v_1,\ldots,v_n$, subtended by the chord 
$\bar{e} = \overline{v_0v_n}$ (see Figure 2), and $\kappa_{H,e}(\pi)$ denotes 
its respective Haantjes curvature. 
%
\begin{figure} 
	\centering
	\includegraphics[height=4cm]{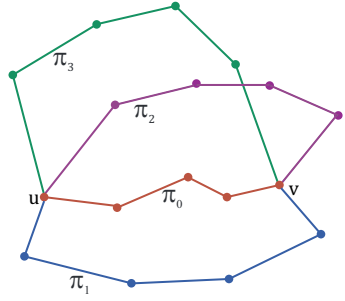}
	\caption{Haantjes-Ricci curvature in the direction $\overline{uv}$,  
		is defined as ${\rm Ric}_H(\overline{uv}) = \sum_1^mK^i_{H,\overline{uv}} = \sum_1^m\kappa_{H,\overline{uv}}(\pi_i)$\, 
		where $\kappa_{H,\overline{uv}}(\pi_i) = \frac{l(\pi_i) - l(\pi_{0})}{ l(\pi_{0})^3}\, i = 1,2,3$; $\pi_0$ being 
		the shortest path connecting the vertices $u$ and $v$.}
\end{figure}

Note that the definition above is much more general than the one based on Menger curvature. 
Indeed, not only is it applicable to cells whose boundary has (combinatorial) length greater 
than three (i.e. not just to triangles), it also does not presume any convexity condition for 
the cells, even in the case when they are realized in some model space, e.g. in $\mathbb{R}^3$.
However, for simplicial complexes endowed with the combinatorial metric, the two notions 
coincide up to a constant. More precisely, in this case, for any triangle $T$,
$\kappa_M(T) \slash \kappa_H(T) = \sqrt{3}/3$. In fact, for
the case of smooth, planar curves Menger and unnormalized
Haantjes curvature coincide in the limit and, furthermore,
they agree with the classical concept. (However, for networks
there is no proper notion of convergence, a fact which allowed
us to discard the  factor 24 in the original definition of Haantjes curvature.)

We can now define, in  analogy with (\ref{eq:Ricci-Menger}), the {\it
	Haantjes-Ricci curvature of an edge} $e$  as 
\begin{equation} \label{eq:Ricci-Haantjes}
{\rm Ric}_H(e) = \sum_{\mathfrak{c} \sim e}K_{H,e}(\mathfrak{c}) =  
\sum_{\mathfrak{c} \sim e}(2\pi -\kappa_{H,e}(\pi))\,;
\end{equation}
where the sum is taken over all the 2-cells $\mathfrak{c}$ adjacent to $e$. 
See Figure 3 for examples of computation of Haantjes-Ricci curvature in networks, in the directed case. 
See Table 1 for the comparison on a number of (undirected) standard planar and spatial grids of the various 
types of Ricci curvature at our disposal.
%
\begin{figure} 
	\centering
	\includegraphics[height=3.4cm]{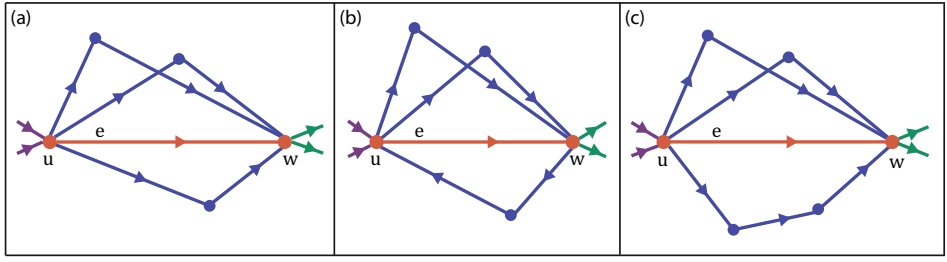}
	\caption{Haantjes-Ricci curvature  ${\rm Ric}_H$ for directed networks, endowed with combinatorial weights: 
		{\bf (a)} ${\rm Ric}_H = 6\pi - 3$, {\bf (b)} ${\rm Ric}_H = 2\pi - 1$, and {\bf (c)} ${\rm Ric}_H = 6\pi - 2 - \sqrt{2}$. 
		Note that while Forman-Ricci curvature is a counter of triangles in simplicial complexes, Haantjes-Ricci curvature 
		represents a counter of all $n$-gones, since each $n$-gone contributes a $\sqrt{n}$ term.}
\end{figure}
%
%
\begin{table}[htb]
	\large
	{
		\begin{center}
			\begin{tabular}{||c||c|c|c|c||} \hline\hline
				{\em Curvature } & {\em Triangular } & {\it Square} & {\it Hexagonal} & {\em Euclidean}\\
				{\em Type} & {\em  Tessellation} & {\it Tessellation} & {\it Tessellation} & {\em Cubulation}
				\\ \hline\hline
				${\rm Ric}_H(e)$ & {4$\pi$ - 2} & {4$\pi$ - 4} & {4$\pi - 2\sqrt{2}$} & {8$\pi - 4\sqrt{2}$} 
				\\ \hline
				${\rm Ric}_{F,r}(e)$ & {-8} & {-2} & {-2} & {-4} 
				\\ \hline
				${\rm Ric}_{F}(e)$ & {-2} & {0} & {4} & {4} 
				\\ \hline
				${\rm Ric}_{O}(e)$ & {1} & {-1} & {--} & {-$\frac{4}{3}$} 
				
				\\ \hline\hline
			\end{tabular}
		\end{center}
	}
	
	\caption{Comparison of undirected curvatures for a number of standard grids (tesselations) of the 
		Euclidean plane and space.}
	\label{table:comparisons}
\end{table}

\subsection{The case of general weights}

We now return to the general case. First let us note that for general
weighted graphs, it is not reasonable to attach area 1 to every
2-cell. However, as discussed  in \cite{HJ,SW19}, it is possible to endow
cells in an abstract weighted network with weights that are both derived from 
the original ones and have a geometric content. 
For instance, in the case of social or biological networks, endowed with the combinatorial metric, 
one can designate to each face, instead of the canonical combinatorial weight equal to 1, 
a weight that ``penalizes" the faces with more edges, thus reflecting the weaker mutual 
connections between the vertices of such a face. 
Thus, it is possible to derive a proper local Gauss-Bonnet formula for such general networks, 
i.e. in a manner that still retains the given data, yet captures the geometric meaning of area, 
volume, etc.

Thus, when considering any such geometric weight $w_g(\mathfrak{c}^2)$ of a 2-cell $\mathfrak{c}^2$, 
the appropriate form of the first term on the left side of (\ref{eq:SmoothGB}) becomes
\[
Kw_g(\mathfrak{c}^2)\,,
\]
and the fitting form of (\ref{eq:K=-kH}) is
\begin{equation} \label{eq:K=-kH+}
K_{H,e}(\mathfrak{c}) = -\frac{1}{w_g(\mathfrak{c}^2)}\left(2\pi - \kappa_{H,e}(\pi)\right)\,;
\end{equation}
Before passing to the problem of extending the above definition to the case of general  
weights, let us note that the observations above regarding Menger curvature for directed networks 
apply also to Haantjes curvature, after properly extending the notion of directed 1-cycles of any 
length and not just of directed triangles (see \cite{SW19} and Figure 1). Again, as for the 
Menger curvature, considering directed networks actually simplifies the problem, in the sense 
that it allows for variable curvature (and not just a constant sign one). 
For  general edge weights, we have  the problem that the {\em total weight} $w(\pi)$ of a path 
$\pi = v_0v_1\ldots v_n$ is not necessarily smaller than the weight of its subtending chord $v_0v_n$ 
\footnote{We suggest the name {\em strong local metrics} for those sets of positive weights that 
	satisfy the generalized triangle inequality $w(v_0v_1\ldots v_nv_{n+1}) < w(v_0v_n)$, for any 
	elementary 1-cycle $v_0v_1\ldots v_nv_0$\,.}, thus Haantjes' definition cannot be applied. 
However, we can turn this  to our own advantage by reversing the roles of $w(\pi)$ and $w(v_0v_n)$ 
in the definition of the Haantjes curvature and assigning a minus sign to  the curvature of cycles 
for which this occurs.

Thus, this approach actually allows us to define a variable sign Haantjes curvature of cycles 
(hence a Ricci curvature as well), even if the given network is not a naturally directed one.

Note that the case when $w(v_0v_1\ldots v_n) = w(v_0v_n)$, i.e. that of zero curvature of the 
2-cell $\mathfrak{c}$ with $\partial \mathfrak{c} = v_0v_1\ldots v_nv_0$ straightforwardly corresponds 
to the splitting case for the path  metric induced by the weights $w(v_iv_{i+1})$.

Indeed, the method suggested above reduces to the use of the path metric, in most of the cases. 
However, both the lack of complete generality of the path metric approach and the sign assignment 
advantage exploited in the preceding paragraph, induced us to prefer the straightforward approach above. 
Of course, one can always pass to the path metric and apply to it the Haantjes curvature. Beyond the 
complications that this might induce in certain cases, it is, in our view, less general, at least from 
a theoretical viewpoint, since it necessitates the passage to a metric.  
However, in the case of most general  weights, i.e. both vertex and edge weights, one has to pass to a metric. 
We find the {\it path degree metric} (see e.g. \cite{Ke}) especially alluring, given that it 
combines simplicity with the capacity of capturing in the discrete context essential geometric 
properties of Riemannian metrics. (However, see also \cite{SA05} for an ad hoc metric devised precisely 
for use on graphs in tandem with Haantjes curvature.)


\subsection{A Further Generalization}

Formulas  (\ref{eq:haantjes2}) and (\ref{eq:haantjes-path}) are
meaningful not only for  a single edge. We can consider any two
vertices $u, v$ that can be connected by a path. Among the  (simple) paths $\pi_1,\ldots, \pi_m$ connecting them, 
the shortest one, i.e. the one for which $l(\pi_{i_0}) = \min\{l(\pi_1),\ldots, l(\pi_m)\}$ is attended represents 
the {\it metric segment} of ends $u$ and $v$. Therefore, given any two such vertices, we can define the 
Haantjes-Ricci curvature in the direction $\overline{uv}$ to be 
\begin{equation}
{\rm Ric}_H(\overline{uv}) =  \sum_1^mK^i_{H,\overline{uv}} = \sum_1^m\kappa_{H,\overline{uv}} (\pi_i)\,
\end{equation}
where $K^i_{H,\overline{uv}}$ denotes the Haantjes-Ricci curvature of the cell $\mathfrak{c}_i$, 
where $\partial \mathfrak{c}_i = \pi_i\pi_0^{-1}$, relative to the direction $\overline{uv}$, and where
\begin{equation}
\kappa_H(\pi_i) = \frac{l(\pi_i) - l(\pi_{0})}{ l(\pi_{0})^3}\,,
\end{equation}
and where the paths $l(\pi_1),\ldots, l(\pi_m)$ satisfy the condition that $\pi_i\pi_0^{-1}$ is an {\it elementary cycle}. 
(This  represents a locality condition in the network setting.)

We conclude this section by noting that both the directed  version of curvature and the one for general weights can be 
extended, {\it mutatis mutandis}, to this generalized definition.

\subsection{Simplified Versions for Simplicial Complexes}

For networks in general, but especially in the case of simplicial complexes, it is useful to notice that 
Formula (\ref{eq:haantjes-path}) for a triangle $T = T(uvw)$ reduces to 
\begin{equation}
\kappa^2_H(T) = \frac{d(u,v) + d(v,w) - d(u,w)}{(d(u,w))^3}\,.
\end{equation}

Haantjes curvature of triangles is thus closely related to two other measures, namely the {\it excess} 
${\rm exc}(T)$ and {\it aspect ratio} ${\rm ar}(T)$, that are defined as follows:
\begin{equation}
{\rm exc}(T) = \max_{v \in \{u,v,w\}}{(d(u,v) + d(v,w) - d(u,w))}\,;
\end{equation}
\begin{equation}
{\rm ar}(T) = \frac{e(T)}{d(T)}\,; 
\end{equation}
where $d(T)$ denotes the diameter of a triangle $T = T(uvw)$.

There are strong connections between the excess, aspect ratio, and curvature. In particular, for the normalized 
Haantjes curvature  introduced above, we have the following relation between the three notions:
\begin{equation}
\kappa_{H}^2(T(v)) = \frac{e(T(v))}{d^3}\,,
\end{equation}
that is 
\begin{equation}
\kappa_{H}(T(v)) = \frac{\sqrt{{\rm ar}(T(v))}}{d}\,.
\end{equation}
Since the factor $\frac{1}{d}$ has the role of ensuring that, in the limit, the curvature of a triangle will 
have the dimensionality of the curvature at a point of a planar curve, the aspect ratio can be viewed as a 
(skewed), un-normalized version of curvature (and Haantjes curvature can be viewed as a {\it scaled} version of 
excess). Thus, since the notion of scale is not of true import in many aspects of network understanding, 
curvature can be replaced by these surrogates. 

Also, for the global understanding of the  shape of networks, it is
useful to compute, as is common in the manifold context, the {\it maximal} excess and {\it minimal} aspect ratio 
over all triangles in the network. 


\section{Conclusions and Further Work}

In the present paper we have shown that, based on two classical
notions of metric curvature, namely the Menger and Haantjes
curvatures, it is easy to define expressive notions of curvature, for
networks and their higher-dimensional generalizations, as well as for
simplicial and clique complexes. In particular, it is possible to
define metric Ricci curvatures for quite general networks -- both
with vertex and edge weights. Furthermore, due to the simple definitions
of the metric curvatures residing at the base of these definition, the
new definitions are computationally efficient (especially those based
on Menger curvature), while being, at the same time, extremely
versatile (in particular those derived from Haantjes curvature). In
fact, for combinatorial polyhedral complexes, it proves to be more
expressive than the (full) Forman curvature, since it takes general
$n$-gones into account, and not just of triangles.

The  metric definition based on the so called {\it Wald metric
	curvature} proposed in \cite{GS} allows for the easy
derivation of convergence results as well as the proof of
theoretical results, such as a polyhedral version of the classical Bonnet-Myers Theorem, 
but it is computationally extremely expensive, rendering it practically 
prohibitive as far as concrete calculations are concerned. 
This is in stark contrast with the simplicity and efficiency of the 
methods developed in the present article. 

Since the permitted length of this article is limited, we could only
introduce  the main ideas and definitions, and could treat neither  possible 
applications, nor deeper theoretical aspects. 
Concerning the implementation aspect is concerned, we see two
immediate and necessary tasks that we hope to develop further. The
first one is a statistically significant comparison of the various
extant notions of Ricci curvature -- Forman, reduced Forman, Ollivier,
Stone, etc.  on large empirical  and model networks. The second one is
the exploration of the clustering and community detection capabilities 
of the Ricci curvature notions introduce herein, and their comparison, 
with the results in \cite{SA05} and \cite{NLLG,SJB} respectively. 
Furthermore, it would be interesting to explore the correlation between 
the notions of curvature introduced herein and hyperbolic embeddings of 
networks. In particular, one would like to explore to what extent the 
curvatures predict values using the inferred hyperbolic distances among nodes 
(points) in embeddings like the one considered in \cite{BPK}, and see how 
much these values agree with or deviate from the curvatures measured 
on the observed network. On the theoretical end of the spectrum, one 
would naturally like to prove analogues of such results as the Bonnet-Myers 
Theorem already mentioned above and, most importantly, of a fitting analogue 
of the fundamental global Gauss-Bonnet Theorem, with  important applications 
in the study of long time evolution of networks \cite{WSJ18}.

%
%

\end{document}